%% file: paper.tex
\documentclass[12pt]{amsart}
\font\emailfont=cmtt10

\usepackage{amsmath,amssymb,amsthm,amsfonts,amscd,flafter,epsf, epsfig,
graphicx, graphs}

\input{basmac}

\input{macinv}

\input{macf}

\input{macknots}
\input{macnew}

\title[{Planar open books and Floer homology}] 
{Planar open books and Floer homology}

\author[Peter Ozsv{\'a}th]{Peter Ozsv\'ath}
\address{Department of
Mathematics, University of California at Berkeley, California 94720 \newline
\indent{\emailfont{petero@math.berkeley.edu}}}
\thanks{PSO was supported by NSF grant number DMS 0234311}

\author[Andr{\'a}s Stipsicz]{Andr{\'a}s Stipsicz}
\address{Department of
Mathematics, Institute for Advanced Study, Princeton, New Jersey 08544 \newline
\indent{\emailfont{stipsicz@math.ias.edu}}}
\thanks{AS was supported by OTKA T49449}

\author[Zolt{\'a}n Szab{\'o}]{Zolt{\'a}n Szab{\'o}} 
\address{Department of
Mathematics, Princeton University, New Jersey 08540 \newline
\indent{\emailfont{szabo@math.princeton.edu}}}
\thanks{ZSz was supported by NSF grant number DMS 0107792}

\begin{document}

\begin{abstract}  
  Giroux has described a correspondence between open book
  decompositions on a 3--manifold and contact structures. In this paper
  we use Heegaard Floer homology to give restrictions on contact
  structures which correspond to open book decompositions with planar
  pages, generalizing a recent result of Etnyre.
\end{abstract} 

\include{intbooks}

\include{background}
\include{proof}

\commentable{ \bibliographystyle{plain} \bibliography{biblio} }
\end{document}

%% file: basmac.tex
\newcommand\commentable[1]{#1}

\newcommand{\HF}{HF}

\newtheorem{theorem}{Theorem}[section]

\newtheorem{cor}[theorem]{Corollary}

\newtheorem{lemma}[theorem]{Lemma}

\newtheorem{remark}[theorem]{Remark}

\def\endproof{\relax\ifmmode\expandafter\endproofmath\else
  \unskip\nobreak\hfil\penalty50\hskip.75em\hbox{}\nobreak\hfil\bull
  {\parfillskip=0pt \finalhyphendemerits=0 \bigbreak}\fi}
\def\endproofmath$${\eqno\bull$$\bigbreak}
\def\bull{\vbox{\hrule\hbox{\vrule\kern3pt\vbox{\kern6pt}\kern3pt\vrule}\hrule}}

\newcommand{\Q}{\mathbb{Q}}

\newcommand{\Z}{\mathbb{Z}}

\newcommand{\Image}{\mathrm{Im}}

\newcommand{\cm}{\cdot}

\newcommand{\CDisk}{D}

\newcommand{\ModSWfour}{\mathcal{M}}
\newcommand{\ModFlow}{\ModSWfour}

\newcommand{\SpinC}{{\mathrm{Spin}}^c}

\newcommand\abuts\Rightarrow
\newcommand\Sym{\mathrm{Sym}}

%% file: macinv.tex
\newcommand\ModSphere{\ModFlow\left({\mathbb S}\longrightarrow 
\Sym^{g-1}(\Sigma_{1})\times \Sym^2(\Sigma_{2})\right)}
\newcommand\ModSpheres\ModSphere

\newcommand\HFred{\HF_{\rm red}}

\newcommand\HFp{\HFb}

\newcommand\HFa{\widehat{HF}}
\newcommand\HFb{HF^+}

\newcommand\UnparModSp{\widehat \ModSp}
\newcommand\UnparModFlow\UnparModSp
\newcommand\Mod\ModSp

\newcommand\spin{\mathfrak s}

\newcommand{\spinc}{\mathfrak s}

\newcommand{\spinct}{\mathfrak t}

\newcommand\ModMaps{\mathcal M}
\newcommand\ModSp\ModMaps

\newcommand\Fp[1]{F^{+}_{#1}}

%% file: macf.tex
\newcommand\Dual{\mathcal D}
\newcommand\Duality\Dual

%% file: macnew.tex
\newcommand\rot{\mathrm{rot}}

\newcommand\cp{c^+}
\newcommand\RightDehn{t}

%% file: intbooks.tex
\maketitle

\section{Introduction}
 
In~\cite{ThurstonWinkelnkemper}, Thurston and Winkelnkemper showed
that an open book decomposition of a 3--manifold $Y$ gives rise in a
natural way to a contact structure over $Y$, and hence that every
3--manifold admits some contact structure. More recently,
Giroux~\cite{Giroux} obtained fundamental results stating a kind of
converse to the Thurston--Winkelnkemper result, showing in effect that
contact structures are in one-to-one correspondence with certain
concretely describable equivalence classes of open book
decompositions. This result brought about a revolution in contact
geometry with repercussions throughout low-dimensional topology.  More
visibly, it is the inspiration for the recently proved embedding
result of Eliashberg~\cite{Eliashberg} and Etnyre~\cite{Etnyre}
leading --- among other results --- to the proof that nontrivial knots
have Property $P$ \cite{KMpropP}.  In another direction, it forms the
foundations for the ``contact invariant'' in Heegaard Floer
homology~\cite{HolDiskContact}, a tool which has been helpful in the
classification of contact structures over certain
3--manifolds~\cite{LiscaStipsiczII, LiscaStipsiczI}. 

But Giroux's construction also raises a number of questions: what
contact geometric properties are reflected by topological properties
of the open book decomposition? Or, more specifically, what types of
contact structures correspond to open book decompositions whose pages
are planar? For example, in \cite{EtnyrePlanarBooks}, Etnyre shows that
all overtwisted contact structures are compatible with planar open
book decompositions, as are all tight structures on lens spaces
\cite{Schonenberger}.  Note also that the Weinstein conjecture has been
verified for all contact structures which admit planar open books
\cite{ACH}.  

In~\cite{EtnyrePlanarBooks}, Etnyre gives the following constraints on
contact structures compatible with planar open book decompositions.
(Here, $b_2^+(X)$ resp. $b_2^0(X)$ denotes the maximal dimension of
any subspace on which the intersection form is positive definite resp.
identically zero.)

\begin{theorem}  [Etnyre,~\cite{EtnyrePlanarBooks}] \label{t:et} 
If $X$ is a symplectic filling of a contact 3--manifold $(Y, \xi)$
which is compatible with a planar open book decomposition then
$b_2^+(X)=b_2^0(X)=0$, the boundary of $X$ is connected and if $Y$ is
an integral homology sphere then the intersection form $Q_X$ is
diagonalizable over the integers.
\qed
\end{theorem}

The above theorem can be used to show that many contact structures do
not admit planar open book decompositions. For example, according
to~\cite{EliashbergThurston}, a taut foliation on a 3--manifold $Y$
gives rise to a contact structure which admits a symplectic semi--filling by
the 4--manifold $[0,1]\times Y$, and hence it admits no compatible
planar open book decomposition.

The aim of the present article is to prove a result analogous to
Theorem~\ref{t:et}, using Heegaard Floer homology~\cite{HolDisk}.
Recall that $\HFp(Y)$ is a 3--manifold invariant which is a module
over the polynomial algebra $\Z[U]$. Moreover, the decreasing sequence
of submodules $\{U^d\cdot \HFp(Y)\}_{d=0}^{\infty}$ stabilizes for
sufficiently large $d$, c.f. Section~\ref{HolDisk:sec:DefHF}
of~\cite{HolDisk}.

\begin{theorem}
  \label{thm:GenusZeroOpenBooks} Suppose that the contact structure
  $\xi $ on $Y$ is compatible with a planar open book
  decomposition. Then its contact invariant $\cp(\xi)\in \HFp (-Y)$ is
  contained in $U^d\cdot \HFp(-Y)$ for all $d\in {\mathbb {N}}$.
\end{theorem}
\begin{remark}
{\rm 
Recall that $\HFred(-Y)$ is defined as $\HFp (-Y)/\Image (U^d)$ for some
sufficiently large $d$. Theorem~\ref{thm:GenusZeroOpenBooks} then
translates to the statement that the contact invariant of a contact
structure on a 3--manifold $Y$ compatible with a planar open book
decomposition vanishes when regarded as an element of the quotient
group $\HFred(-Y)$.}
\end{remark}

Etnyre's theorem can be seen as a consequence of
Theorem~\ref{thm:GenusZeroOpenBooks} (and its method of proof).  There
are, however, other contact structures which do not admit planar open
books as a result of Theorem~\ref{thm:GenusZeroOpenBooks}.  We list
here some consequences of Theorem~\ref{thm:GenusZeroOpenBooks}.

\begin{cor}\label{c:nontorsion}
Suppose that $\cp (\xi )\neq 0$ and the associated spin$^c$ structure
$\spin (\xi )$ is nontorsion $($that is, $c_1 (\spin (\xi ))$ is not a torsion
class$)$. Then $\xi $ does not support a planar open book decomposition.
\end{cor}



\begin{cor}\label{c:stein}
Suppose that the contact 3--manifold $(Y, \xi )$ with $c_1(\spin (\xi ))=0$
admits a Stein filling $(X, J)$ such that $c_1(X,J)\neq 0$. Then
$\xi$ is not supported by a planar open book decomposition.
\end{cor}

\begin{cor}\label{c:knots}
Suppose that $L\subset (S^3, \xi_{st})$ is a Legendrian knot with
zero Thurston--Bennequin invariant. Then the contact structure
$\xi _L$ given by Legendrian surgery along $L$ is not supported by 
planar open book decomposition.
\end{cor}

The next corollary concerns the 2--plane fields underlying contact
structures over rational homology 3--spheres $Y$. 
Recall that homotopy classes of 2--plane fields $\xi$ over $Y$
are classified by a pair of data, their induced $\SpinC$ structure
$\spinc(\xi)$, and a Hopf invariant $d_3(\xi)$, which takes values in
$\Q$ when $b_1(Y)=0$.

According to a theorem of Kronheimer and
Mrowka~\cite{KMcontact}, there are only finitely many homotopy classes
of 2--plane fields which represent symplectically fillable contact
structures. We have the following refinement for contact structures
which are both symplectically fillable and also compatible with planar
open book decompositions.

\begin{cor}
\label{c:2planefields}
Suppose that $Y$ is a rational homology 3--sphere.  The number of
homotopy classes of 2--plane fields which admit contact structures
which are both symplectically fillable and compatible with planar open
book decompositions is bounded above by the number of elements in
$H_1(Y;\Z)$. More precisely, each $\SpinC$ structure $\spinc$ is
represented by at most one such 2--plane field, and moreover, the Hopf
invariant of the corresponding 2--plane field must coincide with the
``correction term'' $d(-Y,\spinc)$.
\end{cor}

For the last part of the above statement, recall that on a rational
homology 3--sphere $Y$ the Floer homology $\HFp(Y,\spinc)$ admits a
grading by $\Q$.  The correction term referred to in the corollary is
the function $d\colon \SpinC(Y)\longrightarrow \Q$, which measures the
minimal degree of any homogeneous element in $\HFp(Y,\spinc)$ which
lies in the image of $U^n$ for all non-negative integers $n$,
cf.~\cite[Section~\ref{AbsGraded:sec:CorrTerm}]{AbsGraded}.  This
function is analogous to the gauge-theoretic invariant
Fr{\o}yshov~\cite{Froyshov}.

This paper is organized as follows. In Section~\ref{sec:Background},
we review some of the background needed for the proof of
Theorem~\ref{thm:GenusZeroOpenBooks}: Giroux stabilizations and open book
decompositions, Heegaard Floer homology, and the contact invariant. In
Section~\ref{sec:Proof} we turn to the proof of
Theorem~\ref{thm:GenusZeroOpenBooks} and its corollaries. The proof of
Theorem~\ref{thm:GenusZeroOpenBooks} falls naturally into two steps, the first
of which is a result about monoid generators for the mapping class
group of a planar surface, and the second of which is the calculation
of the Heegaard Floer homology groups of a family of model
3--manifolds. With Theorem~\ref{thm:GenusZeroOpenBooks} in hand, 
we derive its corollaries stated above, and also
give an alternate proof of Theorem~\ref{t:et}.

%% file: background.tex
\section{Open book decompositions and Heegaard Floer homology}
\label{sec:Background}

The aim of this section is to review the relevant background needed
for our present purposes. In Subsection~\ref{subsec:GirouxStab} we
describe the notion of ``Giroux stabilization'', which gives the
equivalence relation between open book decompositions inducing the
same contact structure. For more on this, see
\cite{EtnyreNotes, Giroux}.  
In Subsection~\ref{s:ss2.2} we review some basic
facts regarding Heegaard Floer homologies and discuss a special class
of 3--manifolds, called {\em $L$--spaces}, whose Heegaard Floer
homology groups are as simple as possible.  In
Subsection~\ref{subsec:Contact}, we describe the invariant
$\cp(\xi)\in\HFp(-Y)$ associated to a contact structure $\xi$ over $Y$,
defined with the help of Giroux's results.

\subsection{Open books}
\label{subsec:GirouxStab}

Let $\phi$ be an automorphism of an oriented surface $F$ with
nonempty boundary, and suppose that $\phi$ fixes $\partial F$.  We
can form the mapping torus 
$$M_\phi=F\times [0,1]/(\phi(x),0)\sim
(x,1)$$ to obtain a 3--manifold which fibers over the circle, and whose
boundary is $\partial F\times S^1$.  There is a canonically associated
closed 3--manifold $Y$ obtained from $M_\phi$ by attaching solid tori
$\partial F\times \CDisk ^2$ using the identifications suggested by the
notation. The data $(F,\phi)$ is called an {\em open book
decomposition of $Y$}, and $\phi$ is called the monodromy of the open
book. Thus, an open book decomposition of $Y$ gives rise to a link
$L\subset Y$ which is fibered, called the {\em binding} of the open
book decomposition, while the fibers of $M_{\phi}$ are called its {\em
pages}.  An open book decomposition is said to be {\em compatible}
with a contact structure $\xi$ (or the given open book {\em supports}
$\xi$) if there is a contact 1--form $\alpha \in \Omega ^1 (Y)$ such
that $\ker \alpha$ is isotopic to $\xi$, $d\alpha$ is a positive
volume form on each page of the open book decomposition and $\alpha$
evaluates positively on a tangent vector of the binding $L$ generated
by the orientation of $L$ compatible with the orientation of the
pages.

The construction of Thurston and Winkelnkemper
\cite{ThurstonWinkelnkemper} associates a compatible contact structure
to an open book decomposition of $Y$.  Indeed, according to recent
work of Giroux ~\cite{Giroux}, every contact structure is induced by
an open book decomposition in this manner and, in fact, Giroux gives an
explicit criterion for when two open books induce isotopic contact
structures.

Let $\phi$ be an automorphism of $F$ fixing its boundary. A {\em
  Giroux stabilization} $(F',\phi')$ of $(F,\phi)$ is a new
  surface-with-boundary $F'$, equipped with an automorphism $\phi'$
  obtained as follows. Let $F'$ be obtained from $F$ by attaching a
  1--handle, and let $\gamma$ be a curve which runs through the
  1--handle geometrically once. The automorphism $\phi'$ is then
  obtained by extending $\phi$ over $F'$ by the identity map over the
  1--handle, and then composing by a right-handed Dehn twist
  $\RightDehn_\gamma$ along $\gamma$, that is,
$$\phi'=\phi\circ \RightDehn_\gamma.$$

It is not hard to see that Giroux stabilizations leave the
3--manifold and indeed the associated contact structure unchanged.
Giroux's theorem~\cite{Giroux} states that two open book
decompositions of $Y$ are compatible with isotopic contact structures
if and only if they can be connected by a sequence of Giroux
stabilizations/destabilizations.

\subsection{Heegaard Floer homologies}
\label{s:ss2.2}
Let $W$ be a connected, oriented four-manifold with two boundary
components, $\partial W = -Y_1 \cup Y_2$. Then, we denote this
cobordism by $W\colon Y_1 \longrightarrow
Y_2$. Recall~\cite{HolDiskFour} that there is an induced ${\mathbb
{Z}}[U]$--equivariant map 
\[
\Fp{W}\colon \HFp(Y_1)\longrightarrow
\HFp(Y_2)
\]
on Heegaard Floer homology.  Recall also~\cite{HolDisk}
that the Heegaard Floer homology $\HFp(Y)$ of a 3--manifold $Y$
naturally splits into summands indexed by $\SpinC$ structures over
$Y$.  

Heegaard Floer homology groups are hard to determine in general.
One very useful calculational tool is the surgery exact triangle
which relates the Heegaard Floer homology groups of three suitably related
three-manifolds.  More specifically, the three 3--manifolds
$Y_1,Y_2,Y_3$ form a {\em triad} if there is a knot $K\subset Y_1$
such that $Y_2,Y_3$ can be given by integer surgeries along $K$,
and the framing for producing $Y_3$ is one higher than the one
giving rise to $Y_2$. 
Let $W_1\colon Y_1\to Y_2$ denote the cobordism specified by the
original framing of $K$ in $Y_1$, and $W_2:Y_2\to Y_3$ specified
by $(-1)$--framed surgery along a normal circle $C$ to $K$. Finally,
$W_3\colon Y_3\to Y_1$ denotes the cobordism induced by 
$(-1)$--surgery along a normal circle $D$ to $C$. Then the surgery 
exact triangle (Theorem~\ref{HolDiskTwo:thm:GeneralSurgery} of~\cite{HolDiskTwo})
takes the form
\vskip.2cm
\[
\begin{graph}(6,2)
\graphlinecolour{1}\grapharrowtype{2}
\textnode {A}(1,1.5){$\HFp (Y_1)$}
\textnode {B}(5,1.5){$\HFp (Y_2)$}
\textnode {C}(3,0){$\HFp (Y_3)$}
\diredge {A}{B}[\graphlinecolour{0}]
\diredge {B}{C}[\graphlinecolour{0}]
\diredge {C}{A}[\graphlinecolour{0}]
\freetext (3,1.9){$F^+_{W_1}$}
\freetext (4.6,0.6){$F^+_{W_2}$}
\freetext (1.4,0.6){$F^+_{W_3}$}
\end{graph}
\]
\vskip.2cm

An {\em $L$--space} is a rational homology 3--sphere with the property
that the map $U\colon \HFp(Y) \longrightarrow \HFp(Y)$ (and hence
$U^d$ for all $d\in {\mathbb {N}}$) is surjective.  The lens space
$L(p,q)$ is an $L$-space, for example, and the connected sum of two
$L$--spaces is also an $L$--space.  For a more thorough discussion on
$L$--spaces see \cite{NoteLens}.  In particular, in 
\cite[Proposition~\ref{NoteLens:prop:LSpaces}]{NoteLens},
the following result is proved using the surgery triangle:

\begin{lemma}\label{l:trio}
If a triad $(Y_1,Y_2,Y_3)$ of rational homology spheres satisfy that
$Y_1,Y_3$ are $L$--spaces and the cobordism $W_3\colon Y_3\to Y_1$ has
$b_2^+(W_3)=1$ then $Y_2$ is an $L$--space. \qed
\end{lemma}
We will use this principle in proving the following:
\begin{theorem}\label{t:zlspace}
The 3--manifold $Z$ given by the Kirby diagram of Figure~\ref{f:kirby}
with $p_i,q_i\geq 1$
\begin{figure}[ht]
\includegraphics[height=5cm]{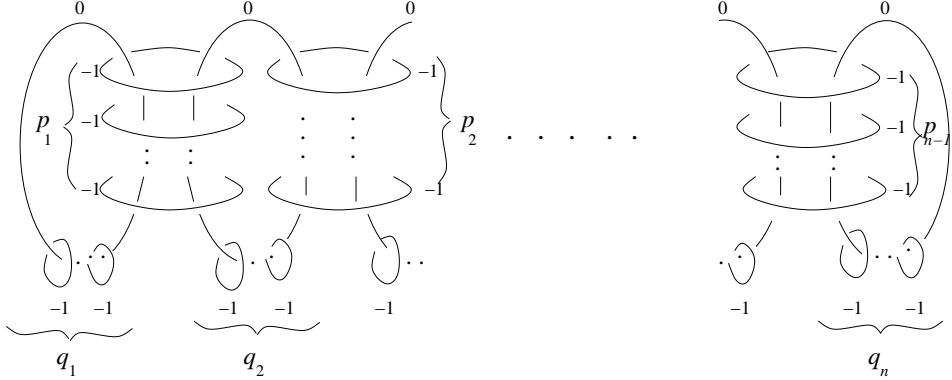}
\caption{\quad Kirby diagram for the 3--manifold $Z$}
\label{f:kirby}
\end{figure}
is an $L$--space.
\end{theorem}
\begin{proof}
We will verify the statement by induction first on the number $n$ of
0--framed unknots and then on $q_n$.  Notice first that in the case of
a single 0--framed knot (by the assumption $q_1\geq 1$) the
3--manifold $Z$ is a lens space, hence the statement easily follows.
Suppose now that the theorem holds for all 3--manifolds of the type
given by Figure~\ref{f:kirby} involving at most $(n-1)$ 0--framed
unknots, and consider $Z$ built with $n$ of those.  
Suppose by induction that for $q_n-1$ the statement is true, and 
consider the triad given by the last $(-1)$--framed unknot $K$ meridional
to the $n^{th}$ 0--framed unknot. If $q_n-1$ is still at least 1, then 
the first element of the triad is an $L$--space by induction on $q_n$,
while the third manifold in the triad (when doing 0--surgery on the
knot $K$) can be given by a surgery diagram of the type given
in Figure~\ref{f:kirby}, now with $(n-1)$ unknots with 0--framing.
Hence our inductive hypothesis shows that it is an $L$--space.
For $q_n=1$ we can easily observe that the first manifold in the
triad (i.e, when we delete the single meridional curve linking the
last 0--framed unknot) admits a presentation of the type of
Figure~\ref{f:kirby}, since the last 0--framed unknot can be canceled againts
one of the $(-1)$--framed circles.
Therefore, in the light of Lemma~\ref{l:trio} the proof of the
theorem follows once we check the condition on the cobordism $W_3$.
In this case it is given by the 2--handle attachment along the dashed
curve $D$ of Figure~\ref{f:cob}. (The solid curves represent 
a Kirby diagram for $Y_3$ in the triad.)
\begin{figure}[ht]
\includegraphics[height=5cm]{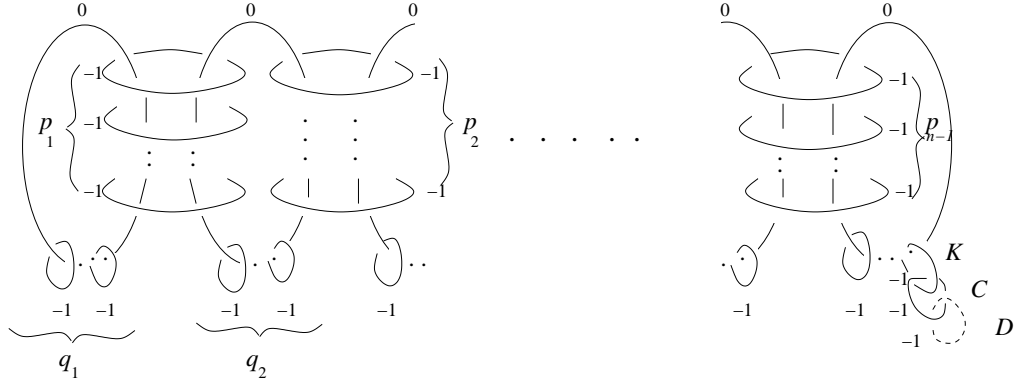}
\caption{\quad Kirby diagram for the cobordism $W_3$}
\label{f:cob}
\end{figure}
Blow down $C$ and slide $D$ over the $n^{th}$ 0--framed unknot, and 
finally cancel $K$ against the $n^{th}$ 0--framed unknot.
In the resulting diagram (shown by Figure~\ref{f:cob2})
the dashed curve represents the cobordism $W_3$.
\begin{figure}[ht]
\includegraphics[height=5cm]{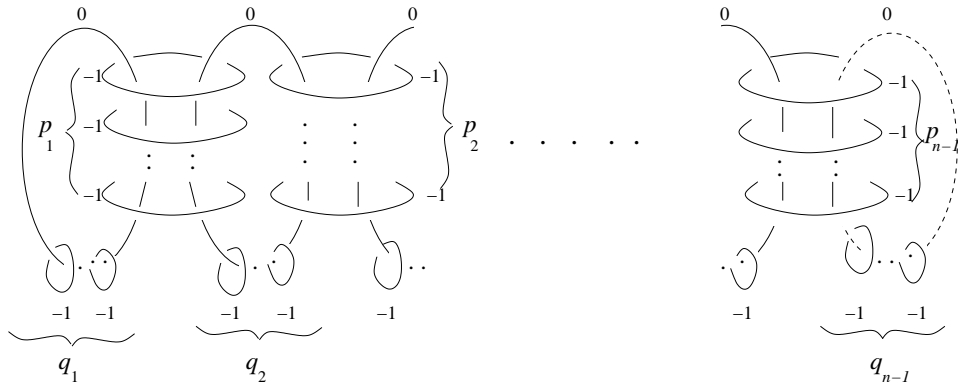}
\caption{\quad Another Kirby diagram for the cobordism $W_3$}
\label{f:cob2}
\end{figure}
Blow down all $(-1)$--circles from the diagram; simple linear algebra
shows that the resulting 4--manifold $X$ will be positive definite. 
Since $W_3$ was shown to admit an embedding into $X$, it follows that
$b_2^+(W_3)=1$, hence the proof of the theorem is complete.
\end{proof}

\begin{remark}
{\rm 
  The surgery diagram of $Z$ determines a planar graph as follows:
  Substitute each 0--framed circle by a vertex $v_i$ ($i=1,\ldots
  ,n$), and connect $v_i$ and $v_{i+1}$ with $p_i$ edges if there are
  $p_i$ $(-1)$--circles linking both the $i^{th}$ and the $(i+1)^{st}$
  unknot.  Finally take an extra vertex $v_{n+1}$ and connect it with
  $q_i$ edges to $v_i$ where $q_i$ is the number of $(-1)$--framed normals
  to the $i^{th}$ 0--framed surgery curve. This planar graph gives
  rise to a connected projection of an alternating link $L$ through
  the black graph of the projection. It can be shown that $Z$ is
  diffeomorphic to the double branched cover $\Sigma (L)$ of $S^3$
  branched along the alternating link $L$. From this observation the
  proof of Theorem~\ref{t:zlspace} is a simple application of
  \cite[Proposition~\ref{BrDCov:prop:AltLink}]{BrDCov}, where it is
  proved that the double branched cover of $S^3$ along a link
  admitting a connected alternating projection is an $L$--space.
}
\end{remark}

\subsection{The contact invariant}
\label{subsec:Contact}
In order to define the contact invariant $\cp (\xi)$, we need one more
observation from Heegaard Floer theory.  Let $Y_0$ be a 3--manifold
which fibers over the circle, with fiber $F$.  Given $i\in\Z$, we can
consider
$$\HFp(Y_0 ,i) = \sum_{\{\spinct\in\SpinC(Y_0)\big| \langle
  c_1(\spinct),[F]\rangle=2i\}}\HFp(Y_0,\spinct).$$
When the genus
$g(F)>1$, then it follows from~\cite{HolDiskContact} that
$\HFp(Y_0,g-1)\cong \Z$, endowed with the trivial action by $\Z[U]$
(i.e., $\HFp(Y_0,g-1)\cong \Z[U]/U\cm \Z[U]$).
Thus, there is a canonical (up to sign) Heegaard Floer homology class
in $\HFp(Y_0)$, which corresponds to a generator of the summand
$\HFp(Y_0,g-1)\subset \HFp(Y_0)$.

Suppose now that $Y$ is a 3--manifold equipped with a contact
structure $\xi$.  We can consider a compatible open book
decomposition. After taking repeated Giroux stabilizations if
necessary, we obtain a new open book decomposition $(F,\phi)$ whose
binding is connected, and whose genus $g(F)$ is greater than one.  By
performing a canonical zero--framed surgery along this connected
binding, we obtain a 3--manifold $Y_0$ which fibers over the circle,
and also a cobordism (obtained by a single 2--handle addition)
$Y\longrightarrow Y_0$. By turning this cobordism upside down, we can
view it as $W\colon -Y_0 \longrightarrow -Y$. The image of a generator
of $\HFp(-Y_0,g-1)\subset \HFp(-Y_0)$ under $\Fp{W}$ in $\HFp(-Y)$ is
the element denoted $\cp(F,\phi)$. It is shown
in~\cite{HolDiskContact} that this element is invariant under Giroux
stabilizations and hence, according to Giroux's theorem, it depends
only on the isotopy class of the underlying contact structure $\xi$
on $Y$; correspondingly, we denote this element by $\cp(\xi)$. (Note
that in~\cite{HolDiskContact} the primary object of study is
a lift of $\cp(\xi)$ from
$\HFp(-Y)$ to $\HFa(-Y)$; we do not need this refinement for our
present applications, however.)

The importance of $\cp (\xi )$ stems from the fact that it seems to
capture interesting contact geometric properties of $\xi$; for
example, $\cp (\xi )=0$ for an overtwisted contact structure, while
$\cp (\xi )\neq 0$ once $\xi $ is Stein fillable. 

The computation of $\cp (\xi )$ can be a very delicate
problem.  A successful scheme of computation rests on the following
result.  To set the stage, let $(F,\phi)$ be a given open book
decomposition for $Y$ with binding $L$, and fix a curve $\gamma\subset
Y-L$ supported in a page of the open book decomposition, which is not
homotopic (in the page) to the boundary. Let $Y_{+1}$ denote the
3--manifold obtained by doing $(+1)$--surgery along $\gamma$ (with
respect to the framing $\gamma$ inherits from the page), and denote
the cobordism defined by this surgery by $X$. Notice that $Y_{+1}$
carries a natural open book decomposition
$(F,\phi\circ\RightDehn_\gamma^{-1})$ and that $-X$ (the 4--manifold
$X$ with its reversed orientation) provides a cobordism $-X\colon -Y
\to -Y_{+1}(\gamma )$.  The following is proved in
\cite[Theorem~\ref{Contact:thm:DehnNaturality}]{HolDiskContact}:

\begin{theorem}
\label{thm:DehnNaturality}
Under the above circumstance, for  the map
\[
\Fp{-X}\colon \HFp(-Y)\longrightarrow \HFp(-Y_{+1}(\gamma))
\]
we have that
\[
\Fp{-X}(\cp(F,\phi))=\pm \cp(F,\phi\circ\RightDehn_\gamma^{-1}).
\]
\qed
\end{theorem}

A contact structure $\xi$ over a 3--manifold $Y$ induces a 
$\SpinC$ structure $\spinc(\xi)$, whose first Chern class $c_1(\xi)$
is the first Chern class of the oriented 2--plane field underlying
the contact structure $\xi$. It is shown in~\cite{HolDiskContact} that
$$\cp(\xi)\in\HFp(-Y,\spinc(\xi))\subset \HFp(-Y).$$
Moreover, since
the maps induced by cobordisms are $\Z[U]$-equivariant, it also
follows that $U\cm \cp(\xi)=0.$ In the case where $c_1 (\xi )$ is
torsion, hence $\HFp (-Y, \spin(\xi))$ admits a ${\mathbb {Q}}$--grading
(defined in~\cite[Section~\ref{HolDiskFour:sec:AbsGrade}]{HolDiskFour},
see also~\cite{AbsGraded}), the element $\cp (\xi )$ is a homogeneous element of
degree $-d_3(\xi )$ in $\HFp (-Y, \spin(\xi))$, where $d_3(\xi )$ is
the Hopf invariant of the oriented 2--plane field underlying $\xi$.
If $(X, \omega )$ is a symplectic filling of $(Y, \xi )$ and
$c_1(\spin (\xi ))$ is torsion then $d_3(\xi )$ can be captured as
\[
\frac{1}{4}(c_1^2(X, \omega )-3\sigma (X_0)-2\chi (X_0)),
\]
where $X_0$ is gotten from $X$ by deleting an open 4--ball
from its interior.


%% file: proof.tex
\section{Proof of Theorem~\ref{thm:GenusZeroOpenBooks}.}
\label{sec:Proof}

In view of the naturality of the contact invariant under left--handed
Dehn twists, the proof of
Theorem~\ref{thm:GenusZeroOpenBooks} breaks into two basic steps:
first, we give a simple set of monoid generators of the mapping class
group of a genus zero surface, so that the contact invariant of any
planar open book decomposition is the image of a Floer homology class
of certain model 3--manifolds, and second, we verify that those model
3--manifolds are $L$-spaces. These two steps are the subjects of the
next two subsections; in Subsection~\ref{subsec:Proofs}, we turn to
the proof of Theorem~\ref{thm:GenusZeroOpenBooks} and its corollaries.

\subsection{Monoid generators for the mapping class group of a planar surface}
\label{subsec:Generators}

\begin{figure}[ht]
\includegraphics[height=5cm]{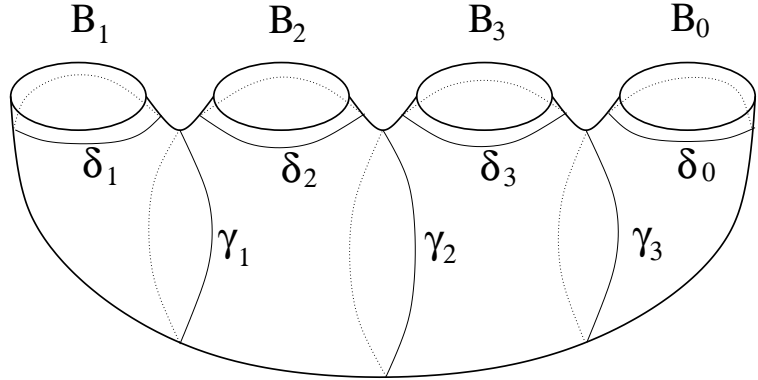}
\caption{\quad The planar surface $S$, in the case where $n=3$.
Here, the automorphisms $\delta_i$ and $\gamma_j$ are
right--handed Dehn twists along the indicated curves.}
\label{f:surface}
\end{figure}

Let $S$ be a compact, planar surface with 
$n+1$ boundary components $B_0,
\ldots , B_n$. Let $\Gamma _S$ denote the mapping
class group of $S$, consisting of diffeomorphisms of $S$ 
which pointwise fix $\partial S$, modulo isotopies, which pointwise
fix $\partial S$.
Let $\delta _i$ denote the
right--handed Dehn twist along a circle parallel to $B_i$ and let
$\gamma _i$ denote the right--handed Dehn twist along a circle
encircling the $i$ boundary components $B_1, \ldots , B_i$, cf. 
Figure~\ref{f:surface}. (Notice that $\delta _1=\gamma _1$ and 
$\delta _0 = \gamma _n$.)

\begin{theorem}
\label{thm:MonoidGenerators}
For an element $x\in \Gamma _S$ there is a decomposition
\[
x=\Pi _{i=1}^n \delta _i ^{n_i}\cdot \Pi _{j=2}^n \gamma _j ^{m_j} \cdot y
\]
where $n_i, m_j$ are positive integers and $y$ can be written as a product
of left--handed Dehn twists.
\end{theorem}
\begin{proof}
It is known that $\Gamma _S$ is generated by Dehn twists. Let
$x\in \Gamma _S$ and write it as a product of Dehn twists 
\[
x=t_1\cdots t_m \in \Gamma _S, 
\]
where $t_i$ denotes a right-- or left--handed Dehn twist along a
simple closed curve in $S$.  Suppose that $t_{i_1}, \ldots , t_{i_k}$
are right--handed Dehn twists.  The idea of the proof is that by using
the ``lantern relation''~\cite{Lantern}
in some related surfaces $S'$ we replace $t_{i_j}$ with
products of left--handed Dehn twists, some other right--handed ones
which have smaller ``complexity'' and with powers of $\delta _i, \gamma
_j$. An inductive argument then provides an expression for $x$
involving left--handed Dehn twists and powers of $\delta _i$'s and
$\gamma _j$'s only.

By complexity we mean the following: suppose that $t_{\alpha}$ is a
Dehn twist along $\alpha \subset S$ encircling the boundary components
$B_{a_1}, \ldots , B_{a_h}$ ($a_1< \ldots < a_h$). Let $c(t_{\alpha
})$ be the maximal element in the difference $\{ 1, 2, \ldots , a_h\}
- \{ a_1, \ldots , a_h\}$.  In particular, if $c(t_{\alpha })=-\infty
$ then $t_{\alpha }= \gamma _{a_h}$.  Among Dehn twists with the same
complexity we say that one is simpler than the other if the maximal
index appearing among the boundary components encircled by it is
smaller than the similar index for the other one. Now using the
lantern relation we can replace any right--handed Dehn twist with
left--handed ones, and with right--handed ones either with smaller
complexity or which are simpler:
\begin{lemma}
Let $t_{\alpha }$ be a given right--handed Dehn twist in $\Gamma _S$.
Then 
\[
t_{\alpha}=t_1\cdot t_2 \cdot \Pi _{i=1}^4 d_i
\]
where $d_i$ are either left--handed or of the form $\delta _j$, $t_1$
has smaller complexity than $t_{\alpha}$ and either
$c(t_2)<c(t_{\alpha })$ or $c(t_2)=c(t_{\alpha })$ and $t_2$ is
simpler than $t_{\alpha }$.
\end{lemma}
\begin{proof}
Suppose that the circle $\alpha $ encircles the boundary
components $B_{a_1}, \ldots , B_{a_h}$, and its complexity
$c(\alpha )\neq -\infty$. Denote the
 circle encircling the same boundary components as $\alpha$
except $B_{a_h}$ by $\beta$.
Define $S'$ as the planar surface we get by substituting the side of
$\beta $ not containing $B_0$ with an annulus.
In  $S'$ we can write down the lantern
relation
\begin{equation} \label{e:lanti}
t_{\alpha }\cdot t_{\lambda _1} \cdot t_{\lambda_2} = 
t_{\beta}\cdot t_{B_{a_h}}\cdot t_{B_{c(\alpha )}}\cdot t_{\lambda _0}
\end{equation}
where we use right--handed Dehn twists everywhere.  (Here $\lambda _1,
\lambda _2$ are the circles encircling $\{ \beta , B_{a_h}\}$ and $\{
B_{c(\alpha )}, B_{a_h}\}$ respectively, while $\lambda _0 $ encircles
$\{ \beta , B_{c(\alpha )}, B_{a_h}\}$.)  Notice that since this
identity holds in $S'$, it will hold in $S$ which can be given by
gluing in an appropriately punctured disk along  the boundary
component of $S'$ corresponding to $\beta$.

After ordering the identity of \eqref{e:lanti} we get an expression
for $t_{\alpha}$ as a product of two left--handed Dehn twists
($t_{\lambda _1}^{-1}$ and $t_{\lambda _2}^{-1}$), two others of the
form $\delta _j$, and $t_{\beta}, t_{\lambda _0}$ . For these latter
two, the complexity of $\lambda _0$ is smaller than that of $\alpha$,
while for $\beta$ we dropped $B_{a_h}$ from the encircled boundary
components, hence either $c(\beta )<c(\alpha)$ or in case $c(\beta
)=c(\alpha)$ then $\beta$ is simpler.  
\end{proof} 

Now induction shows that $x$ can be written as a product of Dehn
twists which are all left--handed except possibly powers of $\delta
_i$ and $\gamma _j$. Since we can conjugate powers of $\delta _i$ and
$\gamma _j$ to the front, and $f^{-1}D_{\alpha}f=D_{f(\alpha )}$ for
any mapping class $f\in \Gamma _S$ and simple closed curve $\alpha$,
we get the desired expression.  In addition, by inserting $\delta
_i\delta _i ^{-1}$ or $\gamma _j \gamma _j ^{-1}$ if necessary, we can
assume that $n_i>0$ and $m_j>0$.  (Recall that $\delta _i, \gamma _j$
are right--handed Dehn twists, hence $\delta _i ^{-1}, \gamma _j ^{-1}$
are left--handed.)
\end{proof}

\subsection{Model 3--manifolds}
\label{subsec:Models}

\begin{theorem}
  \label{thm:ModelManifold}
  Consider a 3--manifold $Z$ which admits a planar open book
  decomposition whose page $S$ has $(n+1)$ boundary components and the
  monodromy map $\phi\in \Gamma _S$ is of the form
  $$\Pi _{i=1}^n \delta _i ^{n_i}\cdot \Pi _{j=2}^n \gamma _j ^{m_j} $$
  with $n_i$ and $m_j$ positive. This 3--manifold 
  is an $L$-space.
\end{theorem}

\begin{proof}
  We give a Kirby calculus description of $Z$. To this end, consider
  first the case when $\phi $ is the identity map. Then the Kirby
  diagram consists of $n$ 0--framed unknots, which can be seen as
  follows. Assume first that $S=D^2$ (i.e., $n=0$), in which case the
  open book decomposition is simply the standard genus--1 Heegaard
  decomposition of $S^3$. 0--surgery on $n$ parallel copies of the
  core circle of the mapping torus $M_1=S^1\times D^2\subset S^3$ then
  provides the 3--manifold corresponding to general $S$ and
  $\phi=1$. Since the page $S$ can be constructed by puncturing a disk
  orthogonal to the surgery curves at the intersections, the curves
  giving rise to the Dehn twists $\delta _i, \gamma _j$ are
  explicitely visible in the picture, cf. Figure~\ref{f:curves}. It is
  known that multiplication of the
\begin{figure}[ht]
\includegraphics[height=5cm]{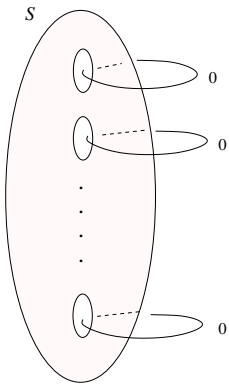}
\caption{\quad Kirby diagram and page for the open book $\phi =1$}
\label{f:curves}
\end{figure}
  monodromy by the right--handed Dehn twists $\delta _i$ or $\gamma
  _j$ changes the 3--manifold by a $(-1)$--surgery along the
  corresponding curve. This observation provides a surgery
  presentation of $Z$. Recall that the boundary components of $S$, and
  so the 0--framed unknots were indexed by $\{ 1, \ldots , n\}$.  To
  get a more convenient presentation of $Z$, slide the $i^{th}$
  0--framed unknot over the $(i+1)^{st}$ (where $i$ runs from 1 to
  $n-1$) in the way that with the natural compatible orientations we
  subtract the corresponding homology classes.  The resulting Kirby
  diagram is now of the type considered in Theorem~\ref{t:zlspace}, where
  it is shown to give rise to an $L$--space, hence the proof is
  complete.
\end{proof}

\begin{remark}
{\rm
Notice that the mapping class $\phi $ also provides a contact 
structure $\xi _{\phi }$ on $Z$. Since $\phi $ is given as
a product of right--handed Dehn twists, it follows 
\cite{Giroux} that $\xi _{\phi }$
is Stein fillable, in particular its contact invariant $\cp (\xi _{\phi })$
is nontrivial in $\HFp (-Z)$. This observation accords with the fact we will
prove in the next Subsection, stating that contact invariants of contact 
structures compatible with planar open book decompositions are images
of $\cp (\xi _{\phi})$'s.
}
\end{remark}

\subsection{Proof of Theorem~\ref{thm:GenusZeroOpenBooks} and 
  its corollaries}
\label{subsec:Proofs} With all the pieces in place, we now turn 
to the proof of Theorem~\ref{thm:GenusZeroOpenBooks}.

\vskip.2cm
\noindent{\bf{Proof of Theorem~\ref{thm:GenusZeroOpenBooks}.}}  Let
$\xi$ be a contact structure on $Y$ specified by a planar open book
decomposition.  Combining Theorem~\ref{thm:MonoidGenerators} with
Theorem~\ref{thm:DehnNaturality}, it follows that there is a
3--manifold $Z$ which admits a planar open book decomposition whose
monodromy $\phi $ has the form
$$\Pi _{i=1}^n \delta _i ^{n_i}\cdot \Pi _{j=2}^n \gamma _j ^{m_j} $$
with $n_i, m_j> 0$ 
and which has the property that $\cp(\xi)$ is the image of the contact
invariant $\cp(\xi_{\phi })\in \HFp (-Z)$ under the map induced by
some cobordism $W\colon -Z\longrightarrow -Y$.
By Theorem~\ref{thm:ModelManifold}, $Z$ is an $L$-space, i.e.
$$U^d\colon \HFp(-Z)\longrightarrow \HFp(-Z)$$ is surjective for all
$d\in {\mathbb {N}}$.  Since maps induced by cobordisms are
$\Z[U]$--module homomorphisms, it follows at once that
$\cp(\xi)\in\HFp(-Y)$ is in the image of
$$U^d\colon \HFp(-Y) \longrightarrow \HFp(-Y),$$ as stated. 
\qed \vskip.2cm




Although Theorem~\ref{t:et} is not a formal consequence of
Theorem~\ref{thm:GenusZeroOpenBooks}, it does follow from the method
of proof, as we illustrate here:

\vskip.2cm
\noindent{\bf{Proof of Theorem~\ref{t:et}.}}  
Suppose that $\xi $ is compatible with the planar open book
decomposition $(S, \phi )$, where $\phi = \Pi _{i=1}^n \delta
^{n_i}\cdot \Pi _{j=2}^n \gamma _j ^{m_j}\cdot y$ with $y$ being the
product of left--handed Dehn twists in $\Gamma _S$. Therefore, the
multiplication of $\phi $ by $y^{-1}$ can be achieved by a sequence of
Legendrian surgeries along $(Y, \xi )$, providing a Stein cobordism
$W$ from $(Y,\xi )$ to some $(Z, \xi ')$, where an open book decomposition
compatible with $\xi '$ is given by $(S, \phi \cdot y^{-1})$. This
observation provides an embedding of the filling $(X, \omega )$ into a
filling of $(Z, \xi ')$. On the other hand, the manifold $Z$ was
proved to be an $L$--space in Theorem~\ref{t:zlspace}, hence by
\cite[Theorem~\ref{Genus:thm:Lspaces}]{GenusBounds} it can be filled
only with symplectic 4--manifolds admitting connected boundary and
vanishing $b_2^+$, implying $b_2^+(X)=b_2^0(X)=0$ and the connectedness
of $\partial X$. 

Suppose now that $Y$ is an integral homology sphere and consider the
cobordism $W\colon Y \to Z$ found above. Notice that by blowing down the
$(-1)$--curves in the diagram of Figure~\ref{f:kirby} we see that 
$Z$ can be given as the bouondary of a positive definite 4--manifold,
therefore $-Z$ can be considered as the boundary of a 
4--manifold $U$ with negative-definite intersection form. 
Then the closed 4--manifold $X\cup _Y \cup W \cup _Z U$ has negative-definite
intersection form $Q_{X\cup_Y W\cup Z}$, and
hence by Donaldson's theorem~\cite{Donaldson}, the form is diagonalizable.
Since $Y$ is an integral homology sphere, we can split the intersection form
$Q_{X\cup_Y W\cup Z}=Q_{X}\oplus Q_{W\cup Z}$. It is an 
easy consequence of a theorem
of Elkies~\cite{Elkies} that $Q_X$ has diagonalizable intersection
form, as well.
This final observation concludes the proof of the theorem.
\qed
\vskip.2cm


\begin{lemma}
Suppose that the contact structure $\xi$ on the rational homology
3--sphere $Y$ is compatible with a planar open book decomposition. If
$X$ is a symplectic filling of $(Y, \xi )$ then $b_1(X)=0$.
\end{lemma}
\begin{proof}
Suppose that $b_1(X)>0$ for a symplectic filling $X$.
If $\vert H_1(Y; {\mathbb {Z}})\vert =n $ then $X$ admits a connected $(n+1)$--fold
unramified cover ${\tilde {X}}$ which is the trivial $(n+1)$--fold cover when
 restricted to $\partial X$. By capping off $n$ of the components of $\partial {\tilde {X}}$
with concave fillings of positive $b_2^+$--invariants, we get a symplectic
filling ${\overline {X}}$ of $(Y, \xi )$ with $b_2^+({\overline {X}})>0$, 
contradicting our previous result.
\end{proof}

\noindent{\bf{Proof of Corollary~\ref{c:nontorsion}.}}
According to
\cite[Proposition~\ref{HFSymp:lemma:AnnihilateHFinf}]{HFSymp}, if
$\spinc$ is a $\SpinC$ structure whose first Chern class is
non--torsion, then there is some integer $m$ with the property that
$U^m \cm \HFp(Y,\spinc)=0$.  On the other hand, according to
Theorem~\ref{thm:GenusZeroOpenBooks}, $\cp(\xi)$ must lie in this
group if $\xi$ supports a planar open book decomposition.
For $c^+(\xi )\neq 0$ this provides a contracidtion, verifying the
corollary
\qed \vskip.2cm

Corollary~\ref{c:stein} follows from
Theorem~\ref{thm:GenusZeroOpenBooks}, together with known properties
of the contact invariants of Stein manifolds (cf.~\cite{LM}
and~\cite{Plamenevskaya}):

\vskip.2cm

\noindent{\bf{Proof of Corollary~\ref{c:stein}.}}
Suppose that $(X, J)$ is a Stein manifold with contact boundary
$(Y, \xi )$, and assume that $c_1(X, J)$ is nonzero in $H^2 (X; {\mathbb {Z}})$
while $c_1(\xi )=0$. Let ${\overline {J}}$ denote the conjugate 
complex structure, inducing ${\overline {\xi }}$ on $Y$. For the induced 
$\SpinC$ structure $\spin (J)$ on $X$ the condition
 $c_1(X, J)\neq 0$ readily implies $\spin (J)\neq
\spin ({\overline {J}})$, hence by \cite{LM} we get that $\xi$ and
${\overline {\xi}}$ are not isotopic, in fact, according to
\cite{Plamenevskaya} we also know that $\cp (\xi )\neq \cp ({\overline
{\xi}})$.  On the other hand, the assumption $c_1(\xi )=0$ implies
that $\xi$ and ${\overline {\xi}}$ induce the same $\SpinC$ structure
$\spin (\xi ) \in \SpinC (Y)$. Since the quotient map
\[ 
R \colon \HFp (-Y, \spin (\xi ))\to HF_{red} (-Y, \spin (\xi ))
\]
has 1--dimensional kernel when restricted to $\ker U$, we get that
at most one of $\cp (\xi )$ and $\cp ({\overline {\xi}})$ can map 
to zero under $R$. Therefore one of $\xi $ and ${\overline {\xi}}$
is not compatible with planar open book decomposition. Since an
open book decomposition of ${\overline {\xi }}$ can be given by
taking an open book decomposition $(F, \phi )$ of $\xi$, 
reversing the orientation
of the page $F$ and inverting $\phi$, we see that 
$\xi $ and ${\overline {\xi }}$ admit planar open book decompositions
at the same time. Since one of them is not compatible with a
planar open book decomposition, the proof is complete.
\qed
\vskip.2cm

\noindent{\bf{Proof of Corollary~\ref{c:knots}.}}
Since tb$(L)=0$, we have that rot$(L)$ is odd, in particular it is not
zero.  Notice that contact $(-1)$--surgery on $L$ also provides a
Stein filling $(X, J)$ for $(S^3_{-1} (L), \xi _L)$ with $c_1(X, J)$
evaluating on the generator of $H_2(X; {\mathbb {Z}})$ as rot$(L)\neq
0$. Since $H_1(S^3 _{-1}(L); {\mathbb {Z}})=0$, it follows that $c_1
(\xi _L)=0$, hence the application of Corollary~\ref{c:stein} implies the
result.
\qed
\vskip.2cm

\noindent{\bf{Proof of Corollary~\ref{c:2planefields}.}}
The contact invariant $\cp(\xi)$ is a homogeneous element of degree
given by $-d_3(\xi)$, with $U \cm\cp(\xi)=0$. 
According to Theorem~\ref{thm:GenusZeroOpenBooks},
the hypothesis that $\xi$ is compatible with a genus zero open book
ensures that
$\cp(\xi)$ lies in $U^m \HFp(-Y,\spinc(\xi))$ for all $m$, while 
the hypothesis that $\xi$ is fillable ensures that 
$\cp(\xi)$ is non-trivial homology class, according to
\cite[Theorem~\ref{Genus:thm:HFaWeakFilling}]{GenusBounds}.
On the other hand, for a rational homology sphere, a  non-trivial homogeneous
homology class in $\HFp(-Y,\spinc)$ which lies simultaneously
in the kernel of $U$ and
the image of $U^m$ is supported in degree $d(-Y,\spinc)$, 
\cite[Section~\ref{AbsGraded:sec:CorrTerm}]{AbsGraded}.
\qed
\vskip.2cm

As an illustration of the above results, let $L$ be a
$(pq-p-q)$--fold stabilization of the Legendrian positive torus knot
$T_{p,q}$ with (maximal) Thurston--Bennequin number $pq-p-q$. 
According to Corollary~\ref{c:knots}, the contact structure $\xi_L$
on $S^3_{-1}(T_{p,q})$ obtained by Legendrian surgery on $L$ is not
compatible with a genus zero open book decomposition.
(In fact, in the cases where $\rot(L)\neq \pm 1$, it is easy to
see that the degree of the contact invariant $-d(\xi_L)$ is different
from zero, which is $-d(S^3_{-1}(T_{p,q}))$; and hence
for those Legendrianizations, one can appeal to the slightly more elementary
Corollary~\ref{c:2planefields}.)

Note that Theorem~\ref{t:et} does not apply to these contact
structures directly, unless someone finds a symplectic filling of
$(S^3_{-1}(L), \xi _L)$ with positive $b_2^+$. We could neither find
nor rule out the existence of such a filling.

%% file: paper.bbl
\begin{thebibliography}{10}

\bibitem{ACH}
C.~Abbas, K.~Cieliebak, and H.~Hofer.
\newblock The {W}einstein conjecture for planar contact structures in dimension
  three.
\newblock math.SG/0409355, 2004.

\bibitem{Donaldson}
S.~K. Donaldson.
\newblock An application of gauge theory to four-dimensional topology.
\newblock {\em J. Differential Geom.}, 18(2):279--315, 1983.

\bibitem{Eliashberg}
Y.~Eliashberg.
\newblock A few remarks about symplectic filling.
\newblock {\em Geom. Topol.}, 8:277--293, 2004.

\bibitem{EliashbergThurston}
Y.~M. Eliashberg and W.~P. Thurston.
\newblock {\em Confoliations}.
\newblock Number~13 in University Lecture Series. American Mathematical
  Society, 1998.

\bibitem{Elkies}
N.~D. Elkies.
\newblock A characterization of the {${Z}\sp n$} lattice.
\newblock {\em Math. Res. Lett.}, 2(3):321--326, 1995.

\bibitem{Etnyre}
J.~B. Etnyre.
\newblock On symplectic fillings.
\newblock {\em Algebr. Geom. Topol.}, 4:73--80, 2004.

\bibitem{EtnyrePlanarBooks}
J.~B. Etnyre.
\newblock Planar open book decompositions and contact structures.
\newblock {\em Internat. Math. Res. Lett.}, 79:4255--4267, 2004.

\bibitem{EtnyreNotes}
J.~B. Etnyre.
\newblock Lectures on open book decompositions and contact structures.
\newblock In {\em Floer homology, gauge theory, and low-dimensional topology},
  2006.

\bibitem{Froyshov}
K.~A. Fr{\o}yshov.
\newblock The {S}eiberg-{W}itten equations and four-manifolds with boundary.
\newblock {\em Math. Res. Lett}, 3:373--390, 1996.

\bibitem{Giroux}
E.~Giroux.
\newblock G\'eom\'etrie de contact: de la dimension trois vers les dimensions
  sup\'erieures.
\newblock In {\em Proceedings of the International Congress of Mathematicians,
  Vol. II (Beijing, 2002)}, pages 405--414, Beijing, 2002. Higher Ed. Press.

\bibitem{Lantern}
D.~L. Johnson.
\newblock Homeomorphisms of a surface which act trivially on homology.
\newblock {\em Proc. Amer. Math. Soc.}, 75:119--125, 1979.

\bibitem{KMcontact}
P.~B. Kronheimer and T.~S. Mrowka.
\newblock Monopoles and contact structures.
\newblock {\em Invent. Math.}, 130(2):209--255, 1997.

\bibitem{KMpropP}
P.~B. Kronheimer and T.~S. Mrowka.
\newblock Witten's conjecture and property {P}.
\newblock {\em Geom. Topol.}, 8:295--310 (electronic), 2004.

\bibitem{LM}
P.~Lisca and G.~Mati{\'c}.
\newblock Tight contact structures and {S}eiberg-{W}itten invariants.
\newblock {\em Invent. Math.}, 129(3):509--525, 1997.

\bibitem{LiscaStipsiczII}
P.~Lisca and A.~I. Stipsicz.
\newblock Ozsv\'ath-{S}zab\'o invariants and tight contact three-manifolds.
  {I}.
\newblock {\em Geom. Topol.}, 8:925--945 (electronic), 2004.

\bibitem{LiscaStipsiczI}
P.~Lisca and A.~I. Stipsicz.
\newblock Seifert fibered contact three-manifolds via surgery.
\newblock {\em Algebr. Geom. Topol.}, 4:199--217 (electronic), 2004.

\bibitem{HFSymp}
P.~Ozsv{\'a}th and Z.~Szab{\'o}.
\newblock Holomorphic triangle invariants and the topology of symplectic
  four-manifolds.
\newblock {\em Duke Math. J.}, 121(1):1--34, 2004.

\bibitem{HolDiskFour}
P.~S. Ozsv{\'a}th and Z.~Szab{\'o}.
\newblock Holomorphic triangles and invariants for smooth four-manifolds.
\newblock math.SG/0110169.

\bibitem{NoteLens}
P.~S. Ozsv{\'a}th and Z.~Szab{\'o}.
\newblock On knot {F}loer homology and lens space surgeries.
\newblock math.GT/0303017, to appear in {\em Topology}.

\bibitem{HolDiskContact}
P.~S. Ozsv{\'a}th and Z.~Szab{\'o}.
\newblock Heegaard {F}loer homologies and contact structures.
\newblock math.SG/0210127, 2002.

\bibitem{AbsGraded}
P.~S. Ozsv{\'a}th and Z.~Szab{\'o}.
\newblock Absolutely graded {F}loer homologies and intersection forms for
  four-manifolds with boundary.
\newblock {\em Advances in Mathematics}, 173(2):179--261, 2003.

\bibitem{BrDCov}
P.~S. Ozsv{\'a}th and Z.~Szab{\'o}.
\newblock On the {H}eegaard {F}loer homology of branched double-covers.
\newblock math.GT/0309170, to appear in {\em Adv. Math.}, 2003.

\bibitem{GenusBounds}
P.~S. Ozsv{\'a}th and Z.~Szab{\'o}.
\newblock Holomorphic disks and genus bounds.
\newblock {\em Geom. Topol.}, 8:311--334 (electronic), 2004.

\bibitem{HolDiskTwo}
P.~S. Ozsv{\'a}th and Z.~Szab{\'o}.
\newblock Holomorphic disks and three-manifold invariants: properties and
  applications.
\newblock {\em Ann. of Math. (2)}, 159(3):1159--1245, 2004.

\bibitem{HolDisk}
P.~S. Ozsv{\'a}th and Z.~Szab{\'o}.
\newblock Holomorphic disks and topological invariants for closed
  three-manifolds.
\newblock {\em Ann. of Math. (2)}, 159(3):1027--1158, 2004.

\bibitem{Plamenevskaya}
O~Plamenevskaya.
\newblock Contact structures with distinct {H}eegaard {F}loer invariants.
\newblock {\em Math. Res. Lett.}, 11(4):547--561, 2004.

\bibitem{Schonenberger}
S.~Schonenberger.
\newblock Private communication.

\bibitem{ThurstonWinkelnkemper}
W.~P. Thurston and H.~E. Winkelnkemper.
\newblock On the existence of contact forms.
\newblock {\em Proc. Amer. Math. Soc.}, 52:345--347, 1975.

\end{thebibliography}
